\documentclass[11pt,a4paper]{amsart}

\usepackage{amsmath,amssymb,amsthm,amsfonts,mathtools}

\usepackage{enumitem, xspace, caption, subcaption, comment,cases}
\usepackage[hidelinks]{hyperref}
\usepackage{stmaryrd}
\usepackage{marvosym}
\usepackage{xcolor}

\newtheorem{theorem}{Theorem}[section]
\newtheorem{lemma}[theorem]{Lemma}

\newtheorem{claim}[theorem]{Claim}
\newtheorem{fact}[theorem]{Fact}

\newtheorem{definition}[theorem]{Definition}

\newtheorem{question}[theorem]{Question}

\theoremstyle{definition}
\newtheorem*{claim*}{Claim}
\theoremstyle{plain}

\makeatletter

\newcommand{\By}[2]{\overset{\mbox{\tiny{#1}}}{#2}}
\newcommand{\ByRef}[2]{   \By{\eqref{#1}}{#2} }

\newcommand{\leBy}[1]{    \By{#1}{\le} }
\newcommand{\geBy}[1]{    \By{#1}{\ge} }

\newcommand{\leByRef}[1]{ \ByRef{#1}{\le} }
\newcommand{\geByRef}[1]{ \ByRef{#1}{\ge} }
\newcommand{\justify}[1]{\fbox{\tiny{#1}}\quad}

\def\endofClaim{\hfill\scalebox{.6}{$\blacksquare$}}
\newcommand{\oldqed}{}
\newenvironment{claimproof}[1][Proof]{
  \renewcommand{\oldqed}{\qedsymbol}
  \renewcommand{\qedsymbol}{\endofClaim}
  \begin{proof}[#1]
}{
  \end{proof}
  \renewcommand{\qedsymbol}{\oldqed}
}

\newcommand{\Mreg}{M^\mathrm{reg}}

\newcommand{\R}{\mathbb{R}}
\newcommand{\N}{\mathbb{N}}
\newcommand{\G}{\mathbb{G}}
\newcommand{\D}{\,\mathrm{d}}

\newcommand{\ca}{\mathcal{A}}
\newcommand{\cb}{\mathcal{B}}
\newcommand{\cc}{\mathcal{C}}
\newcommand{\cd}{\mathcal{D}}
\newcommand{\cE}{\mathcal{E}}

\newcommand{\cp}{\mathcal{P}}
\newcommand{\cq}{\mathcal{Q}}

\newcommand{\cu}{\mathcal{U}}
\newcommand{\cw}{\mathcal{W}}
\newcommand{\cx}{\mathcal{X}}

\newcommand{\cz}{\mathcal{Z}}

\newcommand{\vr}{\mathbf{r}}
\newcommand{\vv}{\mathbf{v}}

\newcommand{\sm}{\setminus}
\newcommand{\eps}{\varepsilon}
\newcommand{\nth}{\varnothing}

\newcommand{\flo}[1]{ \left\lfloor #1 \right\rfloor}
\newcommand{\cei}[1]{ \left\lceil #1 \right\rceil}

\newcommand{\Prob}{\mathbb{P}}
\newcommand{\prob}[1]{\Prob\left( #1 \right)}

\newcommand{\E}{\mathbb{E}}
\newcommand{\Ec}[2]{\E \left(#1\,|\,#2  \right)}

\newcommand{\cutn}[1]{\left\lVert #1\right\rVert_{\square}}
\newcommand{\cutm}{\delta_\square}
\newcommand{\cutnd}{d_\square}
\newcommand{\Lone}[1]{\left\lVert #1\right\rVert_{1}}

\newcommand{\Kernel}{\cw}
\newcommand{\Gra}{{\Kernel_{0}}}

\title{Approximating fractionally isomorphic graphons} 
\author{Jan Hladk\'y}
\author{Eng Keat Hng}
\email{hladky|hng@cs.cas.cz}

\address{Institute of Computer Science of the Czech Academy of Sciences, Pod Vod\'arenskou v\v{e}\v{z}\'i~2, 182~07, Praha~8, Czech Republic}

\thanks{Research supported by Czech Science Foundation Project GX21-21762X and with institutional support RVO:67985807}

\thanks{An extended abstract titled `Fractionally isomorphic graphs and graphons' summarizing this work will appear in the proceedings of Eurocomb2023.}

\begin{document}
\begin{abstract}
Greb\'ik and Rocha [\emph{Fractional Isomorphism of Graphons, Combinatorica 42, pp 365--404 (2022)}] extended the well studied notion of fractional isomorphism of graphs to graphons. We prove that fractionally isomorphic graphons can be approximated in the cut distance by fractionally isomorphic finite graphs. This answers the main question from \emph{ibid}.
As an easy but convenient corollary, we deduce that every regular graphon can be approximated by regular graphs.
\end{abstract}
\maketitle

\section{Introduction} \label{sec:intro}

The theory of limits of dense graph sequences, introduced in~\cite{Lovasz2006,MR2455626} and nicely surveyed in~\cite{LovaszBook}, provides a powerful correspondence between (finite) graphs on the one hand and analytic objects called graphons on the other. In the motivating part of this introduction, we shall refer to this theory rather casually and do not need exact definitions.

The two main results of the theory of dense graph limits say that
\begin{enumerate}[label=(R\arabic{*})]
\item \label{item:graphon-converge-subseq} each sequence of graphs contains a subsequence converging\footnote{convergence is meant here and elsewhere to be in the cut distance (denoted by $\cutm$)} to a graphon, and 
\item \label{item:graphon-seq-converge} each graphon has a sequence of graphs converging to it.
\end{enumerate}
A key question is how various graph properties and parameters relate to their graphon counterparts in the context of~\ref{item:graphon-converge-subseq} and~\ref{item:graphon-seq-converge}. One parameter which is well understood is homomorphism density: by the continuity of homomorphism densities, \emph{every} convergent sequence of graphs has homomorphism densities with respect to each fixed graph which converge to that in the limit graphon. On the other hand, there are properties and parameters of graphons which are more subtle in that they are not reflected in \emph{every} sequence of graphs converging to a certain graphon, but rather in \emph{at least one} sequence. Take for example the concept of vertex-transitivity, which was treated in detail for graphons in~\cite{MR3272377}. It was shown there that the limit graphon of a convergent sequence of vertex-transitive graphs is vertex-transitive; this relates to~\ref{item:graphon-converge-subseq}. Regarding~\ref{item:graphon-seq-converge}, we may ask if every vertex-transitive graphon is the limit of at least one sequence of vertex-transitive graphs; this seems to be a difficult question (see~\cite{PrescribedSymmetries}). Clearly, not every sequence of graphs with a vertex-transitive limit graphon is a sequence of vertex-transitive graphs (hint: delete a single edge). Indeed, many of these more subtle properties concern symmetries.

In this paper we use this point of view to study fractional isomorphism. The concept of fractional isomorphism for graphs was introduced by Tinhofer in 1986~\cite{MR843938}, and several important equivalent but very different-looking definitions were later added by Ramana, Scheinerman, and Ullman~\cite{RamanaScheinermanUllman} and another by Dvořák~\cite{DvorakRecognizeHomomorphism} and independently by Dell, Grohe, and Rattan~\cite{MR3829971}. We refer the reader to Section~6 of~\cite{MR1481157} for a treatise of these concepts and recall here the definition introduced in~\cite{RamanaScheinermanUllman}, which is the most convenient to work with in this paper. An \emph{equitable partition} of a graph $G$ is a partition $\cp = \{P_i\}_{i\in[k]}$ of $V(G)$ into disjoint nonempty parts such that for all $i,j\in[k]$ and $u\in P_i$ we have 
\begin{equation}\label{eq:oks}
\deg_G(u;P_j) = \frac{1}{|P_i|}\sum_{v\in P_i}\deg_G(v;P_j)\;.
\end{equation}
The \emph{parameters} of $\cp$ are a pair $((p_i)_{i\in[k]},(D_{i,j})_{i,j\in[k]})$ where $p_i = |P_i|$ for all $i\in[k]$ and $D_{i,j} = \deg_G(u;P_j)$ for all $i,j\in[k]$ and $u \in P_i$. We say that two graphs $G$ and $H$ are \emph{fractionally isomorphic} if they have equitable partitions $\cp_G$ and $\cp_H$ which can be indexed in a way so that they have the same parameters. As a first example, observe that taking equitable partitions with $k=1$ in the definition yields that any two $d$-regular graphs on $n$ vertices are fractionally isomorphic.

In~\cite{GrebikRocha}, Greb\'ik and Rocha developed a theory of fractional isomorphism for graphons. In particular, they showed that all the equivalent definitions of fractional isomorphism of graphs previously studied have graphon counterparts and that these are indeed all equivalent. It follows from their results that if $\{G_n\}_{n\in\N}$ and $\{H_n\}_{n\in\N}$ are sequences of graphs so that $G_n$ is fractionally isomorphic to $H_n$ (for each $n$), $\{G_n\}_{n\in\N}$ converges to a graphon $W$ and $\{H_n\}_{n\in\N}$ converges to a graphon $U$, then $U$ and $W$ are fractionally isomorphic. Let us sketch a proof of this result. To this end, we first recall another characterization of fractional isomorphism from~\cite{DvorakRecognizeHomomorphism,MR3829971}, which says that two graphs $G$ and $H$ on the same number of vertices are fractionally isomorphic if and only if their homomorphism counts $\hom(T,G)$ and $\hom(T,H)$ coincide for all trees $T$. Hence, if $\{G_n\}_{n\in\N}$, $\{H_n\}_{n\in\N}$, $W$ and $U$ are as above, then by the continuity of homomorphism densities we get that the homomorphism density of every tree $T$ in $W$ and in $U$ are equal. This in turn is a characterization of fractional isomorphism of graphons from~\cite{GrebikRocha}.

As with vertex-transitivity, we cannot reverse this in the strong sense: for two fractionally isomorphic graphons $W$ and $U$ there are many sequences $\{G_n\}_{n\in\N}$ converging to $W$ and $\{H_n\}_{n\in\N}$ converging to $U$ where $G_n$ is not fractionally isomorphic to $H_n$. However, a weak reversal was the main open question in~\cite{GrebikRocha}.

\begin{question}[Question~3.2 in \cite{GrebikRocha}]\label{qu:main}
Let $W$ and $U$ be fractionally isomorphic graphons. Does there exist sequences $\{G_n\}_{n\in\N}$ and $\{H_n\}_{n\in\N}$ of graphs such that $G_n$ is fractionally isomorphic to $H_n$ for each $n\in\N$ and
\[ G_n \to_{\cutm} W \textrm{ and }H_n \to_{\cutm} U\;? \]
\end{question}

The main result of this paper is a positive answer to Question~\ref{qu:main}. In fact, we prove a slightly stronger statement in which we simultaneously approximate an arbitrary (even infinite) family of fractionally isomorphic graphons.

\begin{theorem}\label{thm:main-approx}
Given $\eps>0$ and a family $\cu$ of fractionally isomorphic graphons, there exists $n_0\in\N$ such that for each $n \ge n_0$ there exists a family $\{H_U\}_{U\in\cu}$ of fractionally isomorphic graphs on vertex set $[n]$ such that for each $U\in\cu$ we have $\cutm(U,H_U) \le \eps$.
\end{theorem}

The simplest instance of Theorem~\ref{thm:main-approx} is the regular case. Specifically, suppose that $d\in[0,1]$ and $\cu$ is a family of regular graphons with degree $d$ on a probability space $(\Omega,\pi)$; we say that a graphon $U$ is \emph{regular} with degree $d$ if for almost every $x\in\Omega$ we have $\deg_U(x)=\int_{y \in \Omega} U(x,y) \D\pi=d$. Such graphons are known to be fractionally isomorphic (and it is easy to show this). Hence, our theorem applies and says that there are fractionally isomorphic graphs $\{H_U\}_{U\in\cu}$ of order $n$ which approximate $\{U\}_{U\in\cu}$. This itself does not guarantee that the degree-regularity of the graphs $H_U$. Indeed, for example in our construction we might add an isolated vertex to each graph $H_U$ (which is a step which does not spoil fractional isomorphism). However, it follows from basic properties of the cut distance that most vertices in each graph $H_U$ have degrees around $dn$. Despite this, a closer look at our proof does give that the graphs $H_U$ are in fact $D$-regular for some $D=(1\pm\eps)dn$. This leads to the following theorem.

\begin{theorem}\label{thm:regular-approx}
Suppose that $d\in[0,1]$ and $\cu$ is a family of graphons which are fractionally isomorphic to the constant-$d$ graphon $C_d$. Then for every $\eps>0$ there exists $n_0\in\N$ such that for each $n \ge n_0$ there exists $D\in\N$ and a family $\{H_U\}_{U\in\mathcal{U}}$ of $D$-regular graphs on vertex set $[n]$ such that for each $U\in\mathcal{U}$ we have $\cutm(U,H_U) \le \eps$.
\end{theorem}

The remainder of this paper is organised as follows. We begin with a high level sketch of the strategy of our proof of Theorem~\ref{thm:main-approx} in Section~\ref{ssec:proofidea}. In Section~\ref{sec:prel} we introduce  notation, introduce some technical tools and recall important concepts about graphons, fractional isomorphism for graphons and graphon sampling. We also introduce a version of Szemer\'edi's regularity lemma suited to our purposes. In Section~\ref{sec:proof} we provide our main lemmas and apply them to prove Theorem~\ref{thm:main-approx}. We provide proofs of the main lemmas in Section~\ref{ssec:lemmasproof}. Finally, in Section~\ref{sec:regular} we show how to obtain a proof of Theorem~\ref{thm:regular-approx} from the proof of Theorem~\ref{thm:main-approx} given in Section~\ref{sec:proof}.

\subsection{Proof idea} \label{ssec:proofidea}
Let us describe our proof of Theorem~\ref{thm:main-approx} at a high level and give relevant pointers along the way. The key concept of the theory developed in~\cite{GrebikRocha} is a certain quotient object $W/\cc(W)$ for each graphon $W$; a formal definition is given in Section~\ref{ssec:fracisographonprelim}. The importance of this object is that two graphons $U$ and $W$ are fractionally isomorphic if and only if $U/\cc(U)$ and $W/\cc(W)$ are isomorphic. To start with, we give a simplified description mimicking the definition of equitable partitions of graphs. Let $(\Omega,\pi)$ be a probability space. Assume in this simplification that for a graphon $W \colon \Omega^2 \to [0,1]$ the object $W/\cc(W)$ comes from a suitably chosen finite partition $\cq = \{Q_i\}_{i\in[M]}$ (for some $M\in\N$) of $\Omega$ and has the property that for all $i,j\in[M]$ and $x \in Q_i$ we have 
\begin{equation}\label{eq:musimbalit}
\deg_W(x;Q_j) = \frac{1}{\pi(Q_i)}\cdot\int_{y\in Q_i}\deg_W(y;Q_j)\D\pi\;,
\end{equation}
a property clearly inspired by~\eqref{eq:oks}.
Then $W/\cc(W)$ serves as a counterpart to the profile, that is, it records the measures of the sets $Q_i$ into a vector $\vr = (r_i)_{i\in[M]}$ and the densities of pairs $(Q_i,Q_j)$ into a matrix $\mathcal{D}=(d_{i,j})_{i,j\in[M]}$. That means $r_i=\pi(Q_i)$ for all $i\in[M]$ and $d_{i,j}=\frac{1}{\pi(Q_i)\pi(Q_j)}\int_{(x,y)\in Q_i\times Q_j} W(x,y)$ for all $i,j\in[M]$. Now since fractional isomorphism of graphons is characterized by isomorphic quotients, all the graphons in $\cu$ have the same profile $(\vr,\cd)$.

The standard approach to approximating each $U\in\cu$ is to take a $U$-random graph $H_U\sim\G(n,U)$ for $n$ large. Details about $U$-random graphs are given in Section~\ref{ssec:sampling}; for now it is enough to recall that such a random graph has vertex set $[n]$, which corresponds to a set of $n$ i.i.d.\ uniformly random points $x_1,x_2,\ldots,x_n$ from $\Omega$, and each pair $\{i,j\}$ is inserted as an edge independently with probability $U(x_i,x_j)$. In such a setting, basic concentration tells us that with high probability,
\begin{enumerate}[label=(C\arabic{*})]
	\item \label{item:sketch-sampling-cluster} for each $i\in[M]$ the set $X_i:=\{\ell\in V(H_U) \::\: x_\ell\in Q_i\}$ has size $|X_i|\approx r_i n$, and
	\item \label{item:sketch-sampling-degrees} for each $i,j\in [M]$ and each $\ell\in X_i$ we have $\deg_{H_U}(\ell,X_j)\approx d_{i,j} r_j n$.
\end{enumerate}

Observe that if the approximate equalities above were exact, then the graphs $\{H_U\}_{U\in\cu}$ constructed would automatically be fractionally isomorphic with the collection $\{X_i\}_{i\in[M]}$ playing the role of an equitable partition. Hence, we shall modify each $H_U$ slightly by adding $o(n)$ vertices and $o(n^2)$ edges to obtain exact control on the number of vertices in each part and on degrees between pairs of parts. To get exact control on the number of vertices in each part, for each $i\in[M]$ we add a set $Y_i$ of newly created vertices with $|Y_i| = \flo{r_in(1 + \beta)} - |X_i|$ (for some small $\beta>0$); the vertex set of the modified graph $G_U$ shall be $(X_1\sqcup Y_1)\cup (X_2\sqcup Y_2)\sqcup \ldots\sqcup (X_M\sqcup Y_M)$.

To get exact control on the degrees, we would like to add edges to obtain $G_U$ where (for some $\alpha>0$ small)
\begin{enumerate}[label=(G\arabic{*})]
	\item \label{item:sketch-exact-cluster} for each $i,j\in[M]$ with $i\neq j$ and each $\ell\in X_i\cup Y_i$ we have $\deg_{G_U}(\ell,X_j\cup Y_j) = \flo{(1 + \alpha)d_{i,j}r_jn}$, and
	\item \label{item:sketch-exact-degrees} for each $i\in[M]$ the graph $G_U[X_i\cup Y_i]$ is regular with degree $\flo{(1 + \alpha)d_{i,i}r_in}$.
\end{enumerate}
To achieve the exact degrees in~\ref{item:sketch-exact-cluster}, for each $i,j\in[M]$ with $i\neq j$ we shall add a set of edges in the bipartite pair $(X_i\cup Y_i,X_j\cup Y_j)$ which are incident to $Y_i\cup Y_j$; we obtain these edges by applying a theorem of Gale and Ryser. To achieve the exact degrees in~\ref{item:sketch-exact-degrees}, for each $i\in[M]$ we shall add a set $F_{i,i}$ of edges on $X_i\cup Y_i$ which are incident to $Y_i$; some of these edges are obtained in the bipartite pair $(X_i,Y_i)$ by applying the theorem of Gale and Ryser referenced above, while the remainder of these edges are obtained on $Y_i$ by applying a theorem of Erd\H{o}s and Gallai. Finally, note that since the initial graph contained no edges of the types mentioned above, the addition of such edges does not create multiedges.

The key lemmas that describe the construction outlined in the previous paragraphs are Lemmas~\ref{lem:graphon-frac-iso-approx}--\ref{lem:buffer-precise} (which are all already tailored to the general setting, which we address below). Lemma~\ref{lem:graphon-frac-iso-approx} encapsulates a cleaning procedure for graphons whose purpose goes beyond our simplified setting, Lemma~\ref{lem:sampling-outcome} formalizes the outcome of the random graph sampling procedure (represented by~\ref{item:sketch-sampling-cluster} and~\ref{item:sketch-sampling-degrees} in our simplified setting), and Lemma~\ref{lem:buffer-precise} deals with the subtle procedure of adding vertices and balancing degrees.

In the general setting, $U/\cc(U)$ is defined in terms of a suitable sub-sigma-algebra and does not usually correspond to a finite partition. To see that such a setting is unavoidable in the theory of Grebík and Rocha, observe that vertices inside one set $Q_i$ above must have the same overall degree; on the other hand, it is easy to construct graphons such that any particular degree is achieved only on a nullset of $\Omega$.\footnote{For example, take a graphon $W \colon [0,1]^2\rightarrow [0,1]$ defined by $W(x,y)=(x+y)/2$.} Even though the Grebík--Rocha theory is inherently infinitesimal, for the steps described in~\ref{item:sketch-sampling-cluster}, \ref{item:sketch-sampling-degrees}, \ref{item:sketch-exact-cluster} and~\ref{item:sketch-exact-degrees} it is essential that we work with a finite partition. To bridge the gap between the two settings, we employ Szemer\'edi's regularity lemma. Indeed, it is well-known that if $\cq = \{Q_i\}_{i\in[M]}$ is an $\eps$-regular partition of $U$, then in particular we have an approximate version of~\eqref{eq:musimbalit} for most vertices $x\in Q_i$. All in all, we can essentially return to the finite setting as above, with many additional technical modifications to treat various imprecisions that are unavoidable in connection with the regularity lemma.

\section{Preliminaries} \label{sec:prel}

\subsection{Notation}

Write $\N$ for the set of positive integers and $\N_0$ for the set $\N\cup\{0\}$. For $a\in\N_0$ write $[a]$ for the set $\{1,\dots,a\}$ and $[a]_0$ for $[a]\cup\{0\}$. For $\ell\in\N_0$ and a set $S$ write $\binom{S}{\ell}$ for the set of all subsets of $S$ of size $\ell$. For a set $S$ write $S^{(2)}=\{(i,j) \in S^2:i \ne j\}$. For real numbers $x$ and $y$, we write $x\pm y$ to mean an appropriate real number $z$ satisfying $x-y\leq z\leq x+y$.

\subsection{Tools}

In this subsection we collect some tools. We first state a form of the Chernoff bounds.

\begin{theorem}[{\cite[Theorem 2.1]{JansonLuczakRucinskiBook}}] \label{thm:chernoff-JLR}
Suppose that $X$ is a sum of independent and identically distributed Bernoulli random variables. Then for $0 \le \eps \le 3/2$ we have
\[ \prob{X \ge (1 + \eps)\E[X]} \le e^{-\eps^2\E[X]/3} \textrm{ and } \prob{X \le (1 - \eps)\E[X]} \le e^{-\eps^2\E[X]/3}. \]
\end{theorem}

Next, we recall some tools related to degree sequences in graphs. Recall that a sequence $(c_i)_{i\in [n]}$ is \emph{graphic} if there exists a graph whose degree sequence is $(c_i)_{i\in [n]}$. Graphic sequences are characterized by the Erdős--Gallai theorem.

\begin{theorem}[Erd{\H o}s--Gallai] \label{thm:erdos-gallai}
A sequence $(d_i)_{i\in[n]}$ of nonnegative integers in nonincreasing order is graphic if and only if the following hold.
\begin{enumerate}[label=(EG\arabic{*})]
\item\label{itm:EG1} $\sum_{i=1}^{n}d_i$ is even.
\item\label{itm:EG2} For each $k\in[n]$ we have $\sum_{i=1}^{k}d_i \le k(k-1) + \sum_{i=k+1}^{n}\min\{d_i,k\}$.
\end{enumerate}
\end{theorem}

Next, we state the Gale--Ryser theorem~\cite{MR91855,MR87622}, which characterizes the degree sequences that are realizable as bipartite graphs. More precisely, we say that two sequences $(a_i)_{i\in[n]}$ and $(b_i)_{i\in[m]}$ are \emph{bigraphic} if there exists a bipartite graph with parts $A$ and $B$ with $|A|=n$ and $|B|=m$ so that the degree sequence in $A$ matches $(a_i)_{i\in[n]}$ and the degree sequence in $B$ matches $(b_i)_{i\in[m]}$.

\begin{theorem}[Gale--Ryser] \label{thm:gale-ryser}
A pair $(a_i)_{i\in[n]}$ and $(b_i)_{i\in[m]}$ of sequences of nonnegative integers, with the former in nonincreasing order, is bigraphic if and only if the following hold.
\begin{enumerate}[label=(GR\arabic{*})]
\item \label{itm:GR1} We have $\sum_{i=1}^{n}a_i=\sum_{i=1}^{m}b_i$.
\item \label{itm:GR2} For each $k\in[n]$ we have $\sum_{i=1}^{k}a_i \le \sum_{i=1}^{m}\min\{b_i,k\}$.
\end{enumerate}
\end{theorem}

\subsection{Graphons}
Our notation for graphons mostly follows~\cite{LovaszBook}.
Throughout the paper, $(\Omega,\cb)$ is a standard Borel space endowed with a Borel probability measure $\pi$. We emphasize that $\pi$ may contain atoms.
Write $\Kernel(\Omega)$ for the set of all bounded symmetric measurable functions $W \colon \Omega^2 \to \R$ and call the elements of $\Kernel$ \emph{kernels}. \emph{Graphons}, which are denoted by $\Gra(\Omega)$, are then those kernels whose range is contained in $[0,1]$.

Suppose that $\Lambda$ is another probability space (with an implicit sigma-algebra and measure).
Given a graphon $W\in\Gra(\Omega)$ and a measure-preserving map $\phi \colon \Lambda \to \Omega$, we define a graphon $W^{\phi}\in\Gra(\Lambda)$ by $W^{\phi}(x,y) = W(\phi(x),\phi(y))$.
We say that two graphons $W\in\Gra(\Omega)$ and $U\in\Gra(\Lambda)$ are \emph{isomorphic} if there exists a measure-preserving bijection $\phi \colon \Lambda \to \Omega$ such that $W^{\phi} = U$.

Assume that $\Omega$ is an atomless probability space. Given a graph $F$ on vertex set $V$ we can construct its \emph{graphon representation} $W_F\in\Gra(\Omega)$ as follows. We partition $\Omega$ into sets $(\Omega_v)_{v\in V}$ of individual measures $1/|V|$. For each $(u,v)\in V \times V$ set the graphon $W_F$ to be equal to $1$ on $\Omega_u\times\Omega_v$ if $uv \in E(F)$ and $0$ otherwise. Note that $W_F$ depends on the choice of the partition $(\Omega_v)_{v\in V}$ and hence is not unique. However, $W_F$ is unique up to isomorphism.

By identifying kernels which are equal almost everywhere, we may view $\Kernel(\Omega)$ as a subset of various normed spaces, including $L^{\infty}(\Omega^2)$, and consider corresponding notions of distance on $\Kernel(\Omega)$. 
One of the most combinatorially interesting norms on $\Kernel(\Omega)$ is the \emph{cut norm}, which is given by
\begin{equation} \label{eq:defcn}
\cutn{W} = \sup_{S,T \subseteq \Omega} \left| \int_{S\times T}W \D\pi^2 \right|\;,
\end{equation}
where the supremum is taken over all measurable subsets $S,T \subseteq \Omega$. This gives rise to the \emph{cut norm distance} on $\Kernel(\Omega)$, which plays a central role in graph limit theory and is given by
\begin{equation}\label{eq:defcnd}
\cutnd(U,W) = \cutn{U-W}\;.
\end{equation}
The \emph{cut distance} is inspired by the cut norm distance. Suppose that we have two probability measure spaces $\Omega$ and $\Lambda$, and let $W\in\Gra(\Omega)$ and $U\in\Gra(\Lambda)$. Then we set
\begin{equation}\label{eq:defcd}
\cutm(U,W) = \inf_{\phi} \cutnd(U^{\phi},W)\;,
\end{equation}
where the infimum is taken over all measure-preserving bijections $\phi \colon \Omega \to \Lambda$. Using the concept of graphon representation, we can define $\cutm(U,G):=\cutm(U,W_G)$ and $\cutm(H,G):=\cutm(W_H,W_G)$ for graphs $H$ and $G$.

A deletion of a small number of vertices does not significantly alter a graph in the cut distance. This is stated in the next lemma. This lemma is probably well-known but we could not find it in literature.
\begin{lemma}\label{lem:zoom}
Suppose that $G$ is a graph on $n$ vertices and that $U\subset V(G)$ is a set of $m$ vertices. Then we have $\cutm(G,G[U])\le 2(1-\frac{m}n)$.
\end{lemma}
\begin{proof}
We shall compare a graphon representation of $G$ on an atomless probability space $\Omega$ with one of $G[U]$ on the same probability space. Let us consider a partition $\Omega=\bigsqcup_{v\in V(G)} \Omega_v$ into sets of individual measures $\frac1n$ and a partition $\Omega\setminus\bigcup_{v\in U}\Omega_v=\bigsqcup_{w\in U}\Lambda_w$ into sets of individual measures $\frac{n-m}{nm}$. In particular, for each $w\in U$ we have $\pi(\Omega_w\cup\Lambda_w)=\frac{1}{m}$. Now let $W_G$ be the graphon representation of $G$ with respect to the partition $(\Omega_v)_{v\in V(G)}$, and let $W_{G[U]}$ be the graphon representation of $G[U]$ with respect to the partition $(\Omega_v\cup \Lambda_v)_{v\in U}$. In particular, $W_G$ and $W_{G[U]}$ are equal everywhere except possibly on the set $\left(\Omega\times(\bigcup_{w\in U}\Lambda_w)\right) \cup \left((\bigcup_{w\in U}\Lambda_w)\times\Omega\right)$. This set has measure $2(1-\frac{m}n)-(1-\frac{m}n)^2$.
\end{proof}

Let $(\Omega,\cb,\pi)$ be a standard probability space. The \emph{generalized degree function} of a graphon $W$ is the function $\deg_W \colon \Omega \times \cb \to \R$ given by
\begin{equation}\label{eq:generalizeddegree}
\deg_W(x;Y) = \int_{y \in Y} W(x,y) \D\pi\;.
\end{equation}

The last piece of background information regarding the standard theory of graphons concerns their spectral theory. In fact, we need only the basic concept of associating a self-adjoint Hilbert--Schmidt to each graphon, and even that we only need to introduce the concept of invariant sub-sigma-algebras of a graphon. These concepts are used in this paper to a large extent as blackboxes. Each graphon $W\in\Gra(\Omega)$ defines an operator $T_W \colon L^2(\Omega) \to L^2(\Omega)$ by
\[ (T_Wf)(x) = \int_{y\in\Omega} W(x,y)f(y) \D\pi\;, \]
where $f \in L^2(\Omega)$ and $x\in\Omega$. See~\cite[Section 7.5]{LovaszBook} for details.

\subsection{Fractional isomorphism for graphons}\label{ssec:fracisographonprelim}

We recall some concepts related to the theory of fractional isomorphism for graphons developed by Greb\'ik and Rocha in~\cite{GrebikRocha}. Some of these concepts may seem rather cryptic, especially for a reader whose background is not in analysis. That said, a superficial understanding of this theory should suffice to check most of our arguments.

Throughout this section, $(\Omega,\cb)$ is a standard Borel space with a Borel probability measure $\pi$. We say that $\cc\subseteq\cb$ is a \emph{$\pi$-relatively complete sub-sigma-algebra} if it is a sub-sigma-algebra such that we have $\pi(S \triangle T) > 0$ for all $(S,T)\in(\cb\sm\cc,\cc)$. 
Write $\Theta_{\pi}$ for the set of all $\pi$-relatively complete sub-sigma-algebras of $\cb$. 
Given $\cc\in\Theta_{\pi}$, we write $L^2(\Omega,\cc)$ for the set of all functions in $L^2(\Omega)$ which are $\cc$-measurable. 
Let $W\in\Gra(\Omega)$ be a graphon. We say that $\cc\in\Theta_{\pi}$ is \emph{$W$-invariant} if $T_W(L^2(\Omega,\cc)) \subseteq L^2(\Omega,\cc)$. 
Write $\cc(W)$ for the minimum $W$-invariant relatively complete sub-sigma-algebra of $\cb$. The existence of $\cc(W)$ is established by the combination of Definition 5.12 and Proposition 5.13 in~\cite{GrebikRocha}.

Graphons on a standard Borel space have a sufficiently rich structure to admit notions of averaging and quotients. 

\begin{definition} \label{def:condexpec}
Let $W\in\Gra(\Omega)$ and $\cc\in\Theta_{\pi}$ be $W$-invariant. We define $W_{\cc} = \Ec{W}{\cb\times\cc}$.
\end{definition}

We have the following lemma about $W_{\cc}$ from~\cite{GrebikRocha}.

\begin{lemma}[Claim 5.7 in {\cite{GrebikRocha}}] \label{lem:condexpec-asymm-symm}
Suppose that $W\in\Gra(\Omega)$ and $\cc\in\Theta_{\pi}$ is $W$-invariant. Then $W_{\cc} = \Ec{W}{\cb\times\cc} = \Ec{W}{\cc\times\cc}\in\Gra(\Omega)$.
\end{lemma}

The language of quotient spaces gives an alternative formulation for the notion of conditional expectation on graphons. Let us recall Theorem E.1.\ in~\cite{GrebikRocha}. It says that for each $\cc\in\Theta_{\pi}$ there is a standard Borel space $(\Omega/\cc,\cd)$, a Borel probability measure $\pi/\cc$ on $\Omega/\cc$ and a measurable surjection $q_{\cc}\colon\Omega\to\Omega/\cc$ such that $\pi/\cc$ is the pushforward of $\pi$ via $q_{\cc}$. Furthermore, we have a linear isometry $I_{\cc} \colon L^2(\Omega/\cc,\pi/\cc) \to L^2(\Omega,\pi)$ defined by $I_{\cc}f(x) = f(q_{\cc}(x))$ and its adjoint $S_{\cc} \colon L^2(\Omega,\pi) \to L^2(\Omega/\cc,\pi/\cc)$ satisfies $I_{\cc} \circ S_{\cc} = \Ec{-}{\cc}$.

\begin{definition}
Let $W\in\Gra(\Omega)$ and $\cc\in\Theta_{\pi}$ be $W$-invariant. We define $W/\cc = S_{\cc\times\cc}(W_{\cc})$.
\end{definition}

$W/\cc$ is formally defined on $(\Omega\times\Omega)/(\cc\times\cc)$, but it is easy to verify that there is a measure-preserving bijection $i \colon (\Omega\times\Omega)/(\cc\times\cc) \to \Omega/\cc \times \Omega/\cc$ such that $i \circ q_{\cc\times\cc}(x,y) = (q_{\cc}(x),q_{\cc}(y))$ almost everywhere. Hence, we may identify $(\Omega\times\Omega)/(\cc\times\cc)$ with $\Omega/\cc \times \Omega/\cc$ and assume that $W/\cc\in\Gra(\Omega/\cc)$. 

By Lemma~\ref{lem:condexpec-asymm-symm}, we have $I_{\cc\times\cc}(W/\cc)=W_{\cc}$ and
\begin{equation} \label{eq:condexpec-quotient-equal}
W_{\cc}(x,y) = W/\cc(q_{\cc}(x),q_{\cc}(y))
\end{equation}
for almost every $(x,y) \in \Omega\times\Omega$.

Now we state a definition of fractional isomorphism for graphons. It was shown in~\cite[Theorem 1.2]{GrebikRocha} that this is one of several equivalent characterizations of this concept.

\begin{definition}
We say that graphons $W\in\Gra(\Omega)$ and $U\in\Gra(\Lambda)$ are \emph{fractionally isomorphic} if the quotient graphons $W/\cc(W)$ and $U/\cc(U)$ are isomorphic.
\end{definition}

\subsection{Sampling}\label{ssec:sampling}
We recall the notion of an inhomogeneous random graph model $\G(n,W)$ on a given number of vertices $n$ and for an underlying graphon $W\in\Gra(\Omega)$. This procedure generalizes the Erd\H{o}s--R\'enyi binomial random graphs and plays a central role in the theory of graph limits. For details, see~\cite[Chapter~10]{LovaszBook}. Given $n \in \N$ and a graphon $W$, the random graph $\G(n,W)$ is defined as follows. The vertex set of $\G(n,W)$ is $[n]$. To determine the edge set, sample $n$ independent $\pi$-random elements $x_1,\dots,x_n$ from $\Omega$ and connect each pair $\{i,j\}\in\binom{[n]}{2}$ of distinct vertices independently with probability $W(x_i,x_j)$. Note that the combination of this procedure with a partition of $\Omega$ naturally induces a (random) partition of $[n]$ as follows. Given a partition $\cp=\{P_j\}_{j\in J}$ of $\Omega$, the \emph{partition $\cx=\{X_j\}_{j\in J}$ induced by $\cp$} is the partition of $[n]$ in which each cell $X_j$ of $\cx$ is the set of all vertices $i\in[n]$ for which $x_i\in P_j$.

The following lemma tells us that for sufficiently large $n$ the random graph $\G(n,W)$ is close to $W$ with high probability.

\begin{lemma}[{\cite[Lemma 10.16]{LovaszBook}}] \label{lem:sampling}
For each $n\in\N\sm\{1\}$ and each graphon $W$ we have
\[\prob{\cutm(\G(n,W),W) \le \frac{22}{\sqrt{\log_2 n}}} \ge 1 - \exp\left(-\frac{n}{2 \log_2 n}\right)\;. \]
\end{lemma}

\subsection{Regularity}

Our proof involves determining in advance the parameters of the common equitable partition of the fractionally isomorphic graphs we aim to construct. 
A key idea is that we shall obtain the template for this from an application of a regularity lemma to the common isomorphic quotient $W/\cc(W)$. 
Here we give a few definitions and state the regularity lemma we use.

\begin{definition}\label{def:regularpai}
Let $W\in\Gra(\Omega)$ be a graphon and $C,D \subseteq \Omega$ be disjoint measurable sets of positive measure. The \emph{density of $(C,D)$} in $W$ is defined as
\[ d(C,D) = \frac{1}{\pi(C)\pi(D)} \int_{C \times D}W\D\pi^2\;. \]
Let $\delta>0$. We say that $(C,D)$ is \emph{$\delta$-regular} in $W$ if for all measurable subsets $C' \subseteq C$ and $D' \subseteq D$ we have
\begin{equation}
\left| \pi(C')\pi(D')d(C,D) - \int_{C' \times D'}W\D\pi^2 \right| \le \delta\pi(C)\pi(D)\;.
\end{equation}
\end{definition}

The only property of regular pairs we need is the following well-known fact.
\begin{fact}\label{fact:regdeg}
Let $W\in\Gra(\Omega)$ be a graphon, $\delta>0$, $d\in[0,1]$, and $(C,D)$ be a $\delta$-regular pair of density $d$ in $W$. Then for the sets 
\begin{align*}
R^+&:=\left\{ x \in C : \deg_{W}( x ; D ) \le (d-\sqrt{\delta})\pi(D)\right\}\;\mbox{and}\\
R^-&:=\left\{ x \in C : \deg_{W}( x ; D ) \ge (d+\sqrt{\delta})\pi(D)\right\}\;
\end{align*}
we have $\pi(R^+)\le \sqrt{\delta}\pi(C)$ and $\pi(R^-)\le \sqrt{\delta}\pi(C)$.
\end{fact}

Here is a curious question not directly related to the rest of the paper. The only property of regular pairs we need in this paper is that we can control their degrees in the way stated in Fact~\ref{fact:regdeg}. This property, which one might call `degree-homogeneity', is much weaker than the original property of regularity. Could one directly prove a degree-homogeneity analogue of Szemer\'edi's regularity lemma with much better quantitative bounds than those in Szemer\'edi's regularity lemma?

In the setting of the regularity lemma, the edge distribution\footnote{in case of graphs; in the graphon setting, one should perhaps say `value distribution'} is usually controlled only across pairs of clusters, that is, in the spirit of Definition~\ref{def:regularpai}. Here we need to introduce a notion that captures the property of having the same degree\footnote{`having the same degree' is normally called `regularity', but such a notion would lead to a confusion in combination with Definition~\ref{def:regularpai}} \emph{inside} a cluster.

\begin{definition}\label{def:elegance}
Let $W\in\Gra(\Omega)$ be a graphon and $C\subseteq\Omega$ be a set of positive measure. The \emph{density inside $C$} in $W$ is given by
\[ d(C) = \frac1{\pi(C)^2}\int_{C\times C}W \D\pi^2\;. \]
We say that $C$ is \emph{elegant} in $W$ if we have $\deg_W(x;C)=\pi(C)d(C)$ for almost every $x\in C$. Otherwise, we say that $C$ is \emph{inelegant} in $W$. To be consistent with the concept of density of pairs, we sometimes write $d(C,C)$ instead of $d(C)$ and call this quantity the density of $(C,C)$ in $W$.
\end{definition}

\begin{definition}
Let $W\in\Gra(\Omega)$ be a graphon and $\delta>0$. We say that a partition $\cp = \{P_i\}_{i\in[M]}$ of $\Omega$ is a \emph{$\delta$-regular partition} for $W$ if we have
\[ \sum_{(*)}\pi(P_i)^2+\sum_{(**)}\pi(P_i)\pi(P_j)\le \delta\;, \]
where $(*)$ denotes a summation over all $i\in[M]$ for which $P_i$ is inelegant, and $(**)$ denotes a summation over all $(i,j)\in [M]^{(2)}$ for which $(P_i,P_j)$ is not a $\delta$-regular pair in $W$.
\end{definition}

We call the vector $\vv_{\cp}=(\pi(P_i))_{i\in[M]}$ the \emph{cluster footprint of $\cp$}. Given a finite partition $\cp$ of $\Omega$, we shall often consider it ordered with respect to a fixed order, $\cp = \{P_i\}_{i\in[M]}$.

Now we state the version of Szemer\'edi's regularity lemma for graphons which is tailored to our purposes.

\begin{theorem}\label{thm:regularity}
For each $\delta>0$ there exists $\Mreg(\delta)\in\N$ such that for every standard probability space $(\Omega,\cb,\pi)$, every graphon $W\in\Gra(\Omega)$ has a $\delta$-regular partition into at most $\Mreg(\delta)$ parts.
\end{theorem}

Let us highlight two ways in which our version of the regularity lemma differs from the usual statements of the regularity lemma (for graphons). Firstly, our probability space may have atoms. We remark that this is not entirely novel because the traditional regularization of graphs does not split vertices, which can be thought of as atoms. Secondly, aside from controlling the degrees between clusters, we also control the degrees \emph{inside} clusters using the notion of elegance. These two features come together in the sense that the (square of the) measure contained in the atoms of $\Omega$ may be very substantial in our applications, and hence we have to control them. A notable example of a space with large atoms is $\Omega=\Omega'/\cc(W')$ for some regular graphon $W'\in\Gra(\Omega')$. Indeed, $\Omega$ is a singleton with the Dirac measure in this example.

Here we shall sketch the proof of Theorem~\ref{thm:regularity}. First, we enumerate the atoms in $\Omega$ which have individual measures exceeding $\delta/4$ as $A=\{a_1,\ldots,a_m\}$. Then, we can partition $\Omega\setminus A$ into sets $Q_1,\ldots,Q_\ell$ of individual measures between $\delta/4$ and $\delta/2$. Note that each cluster $\{a_i\}$ is automatically elegant. Hence, the sum of squares of measures of inelegant clusters for the partition $\mathcal{Q}_0=A\cup \{Q_1,\ldots,Q_\ell\}$ is at most $\delta/2$, and this stays true for any further refinement. Now we start the usual cluster-refinement procedure with the initial partition $\mathcal{Q}_0$. (Trivial note: atoms do not get split by the refinement procedure.) We repeatedly refine the clusters along irregular pairs until the total contribution of $\pi(C)\pi(D)$, taken over all irregular pairs $(C,D)$, is less than $\delta/2$. The usual index-pumping argument tells us that the complexity of the final partition is bounded. All in all, we have an outcome as needed for Theorem~\ref{thm:regularity}.

\section{Main lemmas and proof of Theorem~\ref{thm:main-approx}} \label{sec:proof}

We begin with two technical definitions for our proof of Theorem~\ref{thm:main-approx}. These definitions are used to control certain key quantities in our proof.

\begin{definition}\label{def:set1}
For constants $\beta,\lambda,\delta,\alpha>0$ we say that we have \emph{$(\beta,\lambda,\delta,\alpha)$-Setup} if
\begin{align}
\label{eq:constants1}
\beta &\le \frac{1}{1000}\;, \\
\label{eq:constants2}
\lambda &\le \frac{\beta^3}{1000}\;, \\
\label{eq:constant3}
\delta &\le \beta^2\lambda^2\;, \textrm{ and}\\
\label{eq:constants4}
\beta-\alpha&\in [20\lambda,\beta^3/5]\;.
\end{align}

For $M,m\in\N$ we say that we have \emph{$(\beta,\lambda,\delta,\alpha,M,m)$-Setup} if in addition to the above we have $m \ge \frac{10^4M^22^{10^3\beta^{-2}}}{\beta^2\lambda^3\delta^4}$.
\end{definition}

\begin{definition}\label{def:set2}
Suppose that $\beta>0$ and $M\in \N$ are given. An $M \times M$ matrix $\cd = (d_{i,j})_{i,j\in[M]}$ is a \emph{$\beta$-robust density matrix} if it is symmetric and for all $i,j\in [M]^2$ we have $d_{i,j}\in [\beta,1-\beta]$ or $d_{i,j}=0$.
\end{definition}
We say that a vector $\mathbf{v}=(\mathbf{v}_i)_{i\in[M]}\in \R^M$ is \emph{unit} if we have $\mathbf{v}_i\in[0,1]$ for all $i\in[M]$ and $\sum_{i\in[M]}\mathbf{v}_i=1$.

Before stating our main lemmas, let us recall the structure of our proof of Theorem~\ref{thm:main-approx}. The proof proceeds in three main steps. 
First, we tidy up each graphon according to a regular partition on the common isomorphic quotient. We set the density to zero in all inelegant clusters and on all pairs which are irregular or have extremely low density. We also impose a uniform density away from $1$ on pairs with density very close to $1$ and smooth over atypical degrees in regular pairs. This is the subject of Lemma~\ref{lem:graphon-frac-iso-approx}.
Next, we apply sampling procedure from Section~\ref{ssec:sampling} to each tidied up graphon with highly controlled degrees to generate graphs which are close approximations and have partitions that resemble common equitable partitions. This is the subject of Lemma~\ref{lem:sampling-outcome}.
Finally, we carefully add vertices and add and/or delete edges so that we exactly attain the parameters of a common equitable partition. This is the subject of Lemma~\ref{lem:buffer-precise} and involves applications of Theorems~\ref{thm:erdos-gallai} and~\ref{thm:gale-ryser}.

\begin{lemma} \label{lem:graphon-frac-iso-approx}
Assume that $(\Omega,\cb,\pi)$ is a standard Borel probability space.
Assume a $(\beta,\lambda,\delta,\alpha)$-Setup. Given a graphon $W'\in\Gra(\Omega)$, there exists an integer $M\le\Mreg(\delta)$, a unit vector $\mathbf{v}\in[0,1]^{M}$ and a $\beta$-robust density matrix $\cd = (d_{i,j})_{i,j\in[M]}$ such that the following hold. Suppose that $W\in\Gra(\Omega)$ is a graphon such that $W/\cc(W)$ is isomorphic to $W'/\cc(W')$. Then there is a graphon $W_{\beta}\in\Gra(\Omega)$ with a partition $\cq = \{Q_i\}_{i\in[M]}$ of $\Omega$ with cluster footprint $\mathbf{v}$ such that the following hold.
\begin{enumerate}[label=(GE\arabic{*})]
\item \label{itm:graphon-frac-iso-deg-conc} For each $(i,j) \in [M]^2$ and each $x \in Q_i$ we have
\[ \deg_{W_{\beta}}( x ; Q_j ) = (1 \pm 4\lambda)d_{i,j}\mathbf{v}_j\;. \]
\item \label{itm:graphon-edit-close} $\cutn{W-W_{\beta}} \le \Lone{W-W_{\beta}} \le 4\beta$.
\end{enumerate}
\end{lemma}

\begin{lemma} \label{lem:sampling-outcome}
Assume that $(\Omega,\cb,\pi)$ is a standard Borel probability space.
Assume a $(\beta,\lambda,\delta,\alpha,M,m)$-Setup. Suppose that we have a $\beta$-robust density matrix $\cd = (d_{i,j})_{i,j\in[M]}$, a unit vector $\mathbf{v}\in[0,1]^{M}$, and a graphon $W\in\Gra(\Omega)$ with a partition $\cq = \{Q_i\}_{i\in[M]}$ of $\Omega$ with cluster footprint  $\mathbf{v}$. Suppose further that for each $(i,j) \in [M]^2$ and each $x \in Q_i$ we have
\begin{equation} \label{eq:W-cell-degree-bound}
\deg_W(x;Q_j) = (1 \pm 4\lambda)d_{i,j}\mathbf{v}_j\;.
\end{equation}
Then there is a graph $F$ on $m$ vertices with a partition $\cx = \{X_i\}_{i\in[M]}$ of $V(F)$ such that the following hold.
\begin{enumerate}[label=(FS\arabic{*})]
\item \label{itm:F-parts-nearly-equal} For each $i\in[M]$ for which $\mathbf{v}_i\ge \frac{\delta}M$ we have $|X_i|=(1\pm\lambda)\mathbf{v}_i m$.
\item \label{itm:F-frac-iso-deg-nearly-equal} For each $(i,j) \in [M]^2$ for which $\mathbf{v}_j\ge \frac{\delta}M$ and each $x \in X_i$ we have 
\[\deg_F(x;X_j) = (1 \pm 6\lambda)d_{i,j}\mathbf{v}_j m\;.\]
\item \label{itm:F-random-close} We have $\cutm(W,W_F) \le \frac{22}{\sqrt{\log_2 m}}$.
\end{enumerate}
\end{lemma}
As we advertised, in our last lemma we want to modify the graph $F$ from Lemma~\ref{lem:sampling-outcome} so that the approximate equalities from~\ref{itm:F-parts-nearly-equal} and~\ref{itm:F-frac-iso-deg-nearly-equal} become exact ones and in particular the bipartite graphs between different parts become biregular. This necessarily runs into issues of integrality. For example, it is not possible to construct a bipartite graph of nontrivial density (i.e.~$\neq0,1$) for parts of coprime number of vertices. To circumvent this, we introduce the following rounding procedures. Let $\Gamma\in\N$. For a real number $x\in\R$, we set $\llbracket x\rrbracket_\Gamma:=\Gamma\cdot\lfloor \frac{x}\Gamma\rfloor$. For a matrix $\ca = (a_{i,j})_{i,j\in[M]}$, its \emph{$\Gamma$-rounding} $\cb = (b_{i,j})_{i,j\in[M]}$ is the matrix in which $b_{i,j}=\lceil \Gamma a_{i,j}\rceil/\Gamma$.
\begin{lemma} \label{lem:buffer-precise}
Assume a $(\beta,\lambda,\delta,\alpha,M,m)$-Setup. Assume that an even number $\Gamma\in\N$ satisfies $\Gamma\ge \delta^{-2}$ and that $m\ge \Gamma M/\delta^2$. Suppose that we have a $\beta$-robust density matrix $\cd = (d_{i,j})_{i,j\in[M]}$, a unit vector $\mathbf{v}\in[0,1]^{M}$, and a graph $F$ on $m$ vertices with a partition $\cx = \{X_i\}_{i\in[M]}$ of $V(F)$ such that~\ref{itm:F-parts-nearly-equal} and~\ref{itm:F-frac-iso-deg-nearly-equal} hold. Let $\cd^* = (d^*_{i,j})_{i,j\in[M]}$ be the $\Gamma$-rounding of the matrix $\frac{1 + \alpha}{1+\beta}\cd$.

Then for each $n\in\N$ satisfying $(1+\beta)m \le n \le (1+2\beta)m$ there is a graph $G$ on $n$ vertices with a partition $\cz = \{Z_i\}_{i\in[M]_0}$ of $V(G)$ such that the following hold.
\begin{enumerate}[label=(GF\arabic{*})]
\item \label{itm:G-parts-equal} For each $i\in[M]$ we have $|Z_i|=\llbracket (1 +\beta)m\mathbf{v}_i\rrbracket_\Gamma$.
\item \label{itm:G-frac-iso-deg-equal} For each $(i,j) \in [M]^2$ for which $\mathbf{v}_i,\mathbf{v}_j\ge \frac{\delta}M$ and each $v \in Z_i$ we have
\begin{equation} \label{eq:G-frac-iso-deg-equal}
\deg_G(v;Z_j) = d^*_{i,j}|Z_j|\;.
\end{equation}
\item\label{itm:G-zero} All the vertices in the cluster $Z_0$ and in each cluster $Z_i$ with $\mathbf{v}_i<\frac{\delta}M$ are isolated.
\item\label{itm:G-F-close} We have $\cutm(W_F,W_G) \le 5\beta$.
\end{enumerate}
\end{lemma}

Now we apply our main lemmas to give the proof of Theorem~\ref{thm:main-approx}.

\begin{proof}[Proof of Theorem~\ref{thm:main-approx}]
Let $\eps>0$. Without loss of generality, we may assume that $\eps \le 10^{-3}$. Set $\beta = \eps/10$, $\lambda = \beta^4$, $\delta = \beta^{10}$ and $\alpha = \beta - 20\lambda$.
Set $n_1 = \cei{\frac{10^4\Mreg(\delta)^2 2^{10^3\beta^{-2}}}{\beta^2\lambda^3\delta^4}}$ and $n_0 = \cei{(1+\beta)n_1}$. We shall show that given $m \ge n_1$, for all $n\in\N$ satisfying $(1+\beta)m \le n \le (1+2\beta)m$ we have a family $\{H_U\}_{U\in\mathcal{U}}$ of fractionally isomorphic graphs on vertex set $[n]$ such that for any $U\in\mathcal{U}$ we have $\cutm(U,H_U) \le \eps$. By taking $m = \flo{\frac{n}{1+\beta}}$, this clearly implies the desired outcome, that is, for all $n \ge n_0$ we have a family $\{H_U\}_{U\in\mathcal{U}}$ of fractionally isomorphic graphs on vertex set $[n]$ such that for any $U\in\mathcal{U}$ we have $\cutm(U,H_U) \le \eps$.

Let $m \ge n_1$. 
Fix an arbitrary $W'\in\cu$. Lemma~\ref{lem:graphon-frac-iso-approx} returns a number $M\le\Mreg(\delta)$, a robust density matrix $\cd = (d_{i,j})_{i,j\in[M]}$ and a unit vector $\mathbf{v}\in [0,1]^M$. 
Observe that we have a $(\beta,\lambda,\delta,\alpha,M,m)$-Setup. 
Let $n\in\N$ satisfy $(1+\beta)m \le n \le (1+2\beta)m$. 
Set $\Gamma = 2\flo{\frac{\delta^2n_1}{2\Mreg(\delta)}}$; clearly, $\Gamma$ is an even positive integer satisfying $\delta^{-2} \le \Gamma \le \frac{\delta^2m}{M}$.
We shall produce graphs $\{G_W\}_{W\in\cu}$ on $n$ vertices. The key desired fact that all these graphs be fractionally isomorphic follows from that each of them has an equitable partition whose parameters are derived from the (unchanging) matrix $\cd$ and vector $\mathbf{v}$ via~\ref{itm:G-parts-equal}, \ref{itm:G-frac-iso-deg-equal} and~\ref{itm:G-zero}.

Let $W\in\cu$. By definition, $W/\cc(W)$ is isomorphic to $W'/\cc(W')$. Hence, by Lemma~\ref{lem:graphon-frac-iso-approx} we have a graphon $W_{\beta}$ with a partition $\cq = \{Q_i\}_{i\in[M]}$ with cluster footprint $\mathbf{v}$ such that~\ref{itm:graphon-frac-iso-deg-conc} and~\ref{itm:graphon-edit-close} hold. By Lemma~\ref{lem:sampling-outcome} applied with $W_{\beta}$, we have a graph $F$ on $m$ vertices with a partition $\cx = \{X_i\}_{i\in[M]}$ of $V(F)$ such that~\ref{itm:F-parts-nearly-equal}--\ref{itm:F-random-close} hold. By Lemma~\ref{lem:buffer-precise} applied with the graph $F$, there is a graph $G_W$ on $n$ vertices with a partition $\cz = \{Z_i\}_{i\in[M]}$ of $V(G_W)$ such that~\ref{itm:G-parts-equal}--\ref{itm:G-F-close} hold.

The graphs in $\{G_W\}_{W\in\cu}$ are pairwise fractionally isomorphic because they all satisfy~\ref{itm:G-parts-equal}, \ref{itm:G-frac-iso-deg-equal} and~\ref{itm:G-zero}. Let $W\in\cu$. By~\ref{itm:graphon-edit-close}, \ref{itm:F-random-close} and~\ref{itm:G-F-close} we have
\[
\cutm(W,W_{G_W}) \le \cutm(W,W_\beta)+\cutm(W_\beta,W_F)+\cutm(W_F,W_{G_W})
\le 10\beta = \eps \;.
\]
This completes the proof.
\end{proof}

It remains to prove Lemmas~\ref{lem:graphon-frac-iso-approx}--\ref{lem:buffer-precise}.

\subsection{Proofs of main lemmas} \label{ssec:lemmasproof}

First we prove Lemma~\ref{lem:graphon-frac-iso-approx}.

\begin{proof}[Proof of Lemma~\ref{lem:graphon-frac-iso-approx}]
Let $W'$ be a graphon. Set $U' = W'/\cc(W')$. In particular, $U'$ is defined on a probability space $X'=\Omega/\cc(W')$ with measure $\pi_{X'}$. We apply Theorem~\ref{thm:regularity} to the graphon $U'$ and obtain a $\delta$-regular partition $\cp' = \{P'_i\}_{i\in[M]}$ of the space $X'$ for some $M\le\Mreg(\delta)$. Let $\mathbf{v}$ be the cluster footprint of $\cp'$.

For each pair $(i,j)\in[M]^2$ write $d'_{i,j}$ for the density of $(P'_i,P'_j)$ in $U'$. Let $S$ be the collection of pairs $(i,j) \in [M]^2$ where either $i \ne j$ and $(P'_i,P'_j)$ is $\delta$-regular in $U'$, or $i=j$ and $P'_i$ is an elegant cluster in $U'$. Let $T$ be the collection of pairs $(i,j) \in S$ where $\beta \le d'_{i,j} \le 1 - \beta$. Let $T'$ be the collection of pairs $(i,j) \in S$ where $d'_{i,j} > 1-\beta$. For each $(i,j) \in [M]^2$ define the quantity
\begin{equation*}
d_{i,j} =
\begin{cases}
d'_{i,j} & \textrm{ if }(i,j) \in T\textrm{,} \\
1-\beta & \textrm{ if }(i,j) \in T'\textrm{,} \\
0 & \textrm{ otherwise}.
\end{cases}
\end{equation*}
Set $\cd = (d_{i,j})_{i,j\in[M]}$ and note that it is a $\beta$-robust density matrix.

Now let $W$ be a graphon such that $U = W/\cc(W)$ is isomorphic to $U'$. In particular, $U$ is defined on a probability space $X=\Omega/\cc(W)$ with measure $\pi_{X}$. Furthermore, we have a measure-preserving quotient map $q_X:\Omega\rightarrow X$. By the fact that $U$ and $U'$ are isomorphic, we have a measure-preserving bijection $\phi:X'\rightarrow X$ so that $U'=U^\phi$.

Observe that the partition $\cp = \{P_i\}_{i\in[M]}$ of $X$ given by $P_i = \phi(P'_i)$ is $\delta$-regular for $U$, and for each $(i,j) \in [M]^2$ the quantity $d'_{i,j}$ is the density of $(P_i,P_j)$ in $U$. Furthermore, for each $(i,j) \in [M]^{(2)}$ the pair $(P_i,P_j)$ is $\delta$-regular in $U$ if and only if the pair $(P'_i,P'_j)$ is $\delta$-regular in $U'$, and for each $i \in [M]$ the cluster $P_i$ is elegant in $U$ if and only if the cluster $P'_i$ is elegant in $U'$. For each $(i,j) \in T$ define the sets
\begin{align*}
R^{-}_{i,j} & = \left\{ x \in P_i : \deg_{U}( x ; P_j ) \le (d_{i,j}-\sqrt{\delta})\pi_X(P_j)\right\}, \\
R^{+}_{i,j} & = \left\{ x \in P_i : \deg_{U}( x ; P_j ) \ge (d_{i,j}+\sqrt{\delta})\pi_X(P_j)\right\}, \\
R_{i,j} & = R^{-}_{i,j} \cup R^{+}_{i,j} \textrm{ and } P_{i,j} = P_i \sm R_{i,j}.
\end{align*}
In the case where $i\neq j$, we apply Fact~\ref{fact:regdeg} to obtain 
\begin{equation}\label{eq:nal}
\pi_X(R^{-}_{i,j}),\pi_X(R^{+}_{i,j}) \le \sqrt{\delta}\pi_X(P_i)\;.	
\end{equation}
Note that the same inequalities also hold in the case where $i=j$ because we have $\pi_X(R^{-}_{i,j}),\pi_X(R^{+}_{i,j})=0$. Define the graphon $U_{\beta}\in\Gra(X)$ by
\begin{equation*}
U_{\beta}(x,y) =
\begin{cases}
U(x,y) & \textrm{if $(x,y) \in P_{i,j} \times P_{j,i}$ with $(i,j) \in T$,}\\
d_{i,j} & \textrm{if $(x,y)\in P_i\times P_j\setminus P_{i,j} \times P_{j,i}$ with $(i,j) \in T$,}\\
d_{i,j} & \textrm{if $(x,y)\in P_i\times P_j$ with $(i,j) \notin T$}.
\end{cases}
\end{equation*}
Define the partition $\cq = \{Q_i\}_{i\in[M]}$ of $\Omega$ by $Q_i = q_X^{-1}(P_i)$. This definition immediately implies that
\begin{equation}\label{eq:Qmeasurable}
\textrm{$Q_i$ is $\cc(W)$-measurable for each $i\in[M]$\;.}
\end{equation}
For each $(i,j) \in T$ define the sets
\begin{align*}
S^{-}_{i,j} & = \left\{ x \in Q_i : \deg_W( x ; Q_j ) \le (d_{i,j}-2\sqrt{\delta})\pi(Q_j)\right\}, \\
S^{+}_{i,j} & = \left\{ x \in Q_i : \deg_W( x ; Q_j ) \ge (d_{i,j}+2\sqrt{\delta})\pi(Q_j)\right\}, \\
S_{i,j} & = S^{-}_{i,j} \cup S^{+}_{i,j} \quad\textrm{ and }\quad Q_{i,j} = Q_i \sm S_{i,j}.
\end{align*}

Define the graphon $W_{\beta}$ on $\Omega$ as follows.
\begin{equation*}
W_{\beta}(x,y) =
\begin{cases}
W(x,y) & \textrm{ if } (x,y) \in Q_{i,j} \times Q_{j,i} \textrm{ with } (i,j) \in T\textrm{,} \\
U_{\beta}(q_X(x),q_X(y)) & \textrm{ otherwise.}
\end{cases}
\end{equation*}

The following claim helps us show that $W_{\beta}$ satisfies~\ref{itm:graphon-frac-iso-deg-conc} and ~\ref{itm:graphon-edit-close}.

\begin{claim} \label{eq:U-W-mismatch-zero}
For each $(i,j) \in T$ the set $S_{i,j} \sm q_X^{-1}(R_{i,j})$ has measure zero.
\end{claim}
\begin{claimproof}
The proof of the claim does not use regularity or elegance but merely properties of the map $q_X$.
Set $A^{-}_{i,j} = S^{-}_{i,j} \sm q_X^{-1}(R^{-}_{i,j})$. We have
\begin{align*}
&\int_{x\in A^{-}_{i,j}} \deg_{U}(q_X(x);P_j)\D\pi \\
\justify{by~\eqref{eq:condexpec-quotient-equal}}&= \int_{x\in A^{-}_{i,j}} \deg_{W_{\cc(W)}}( x ; Q_j ) \D\pi\\
\justify{by~\eqref{eq:Qmeasurable} and Definition~\ref{def:condexpec}}&= \int_{x\in A^{-}_{i,j}} \deg_W( x ; Q_j ) \D\pi\;.
\end{align*}
We shall make use of both ends of this identity. On the one hand, we have
\begin{equation}\label{eq:end1}
\int_{A^{-}_{i,j}} \deg_W( x ; Q_j ) \le (d_{i,j}-2\sqrt{\delta})\pi(A^{-}_{i,j})\pi(Q_j)
\end{equation}
because $A^{-}_{i,j} \subseteq S^{-}_{i,j}$. On the other hand, we have
\begin{equation}\label{eq:end2}
\int_{A^{-}_{i,j}} \deg_{U}(q_X(x);P_j) \ge (d_{i,j}-\sqrt{\delta})\pi(A^{-}_{i,j})\pi_X(P_j)
\end{equation}
because $q_X(A^{-}_{i,j}) \cap R^{-}_{i,j} = \nth$. Since $q_X$ is measure-preserving, we have $\pi(Q_j)=\pi_X(P_j)$. By comparing~\eqref{eq:end1} and~\eqref{eq:end2}, we obtain $\pi(A^{-}_{i,j}) = 0$. Finally, note that an analogous argument for $A^{+}_{i,j} = S^{+}_{i,j} \sm q_X^{-1}(R^{+}_{i,j})$ implies $\pi(A^{+}_{i,j}) = 0$, so the claim follows.
\end{claimproof}

Let us now verify that $W_{\beta}$ satisfies~\ref{itm:graphon-frac-iso-deg-conc}. We shall consider two cases, which we denote by ($\heartsuit$) and ($\spadesuit$).

\underline{Case ($\heartsuit$):} Suppose $(i,j) \in T$ and $x \in Q_{i,j}$. By the definition of $T$ and~\eqref{eq:constant3} we have $d_{i,j}\ge\beta\ge\sqrt{\delta}/\lambda$. Hence, by the definition of $Q_{i,j}$ we have
\begin{equation} \label{eq:degconc}
\deg_W(x;Q_j) = d_{i,j}(1 \pm 2\lambda)\pi(Q_j)\;.
\end{equation}
We also have
\begin{equation}\label{eq:aso}
	\pi(S_{j,i}) \leBy{C\ref{eq:U-W-mismatch-zero}} \pi_X(R_{j,i}) \leByRef{eq:nal} 2\sqrt{\delta} \pi_X(P_j)\le 2d_{i,j}\lambda\pi_X(P_j) = 2d_{i,j}\lambda\pi(Q_j)\;.
\end{equation}
Since by definition $W(x,y)$ and $W_{\beta}(x,y)$ are equal for $y\in Q_{j,i}$, by the definition given in~\eqref{eq:generalizeddegree} we obtain
\[ \left|\deg_{W}(x;Q_j) - \deg_{W_{\beta}}(x;Q_j)\right| = \left| \int_{y \in S_{j,i}} W(x,y) - W_{\beta}(x,y) \D\pi \right| \le \pi(S_{j,i})\;. \]
Therefore, we have
\begin{align*}
\deg_{W_{\beta}}(x;Q_j) &=\deg_W(x;Q_j)\pm \pi(S_{j,i})\\
\justify{\eqref{eq:degconc} and~\eqref{eq:aso}}&= d_{i,j}(1 \pm 4\lambda)\pi(Q_j)\;,
\end{align*}
as required for~\ref{itm:graphon-frac-iso-deg-conc} in this case.

\underline{Case ($\spadesuit$):} Suppose that either $(i,j) \in T$ and $x \in S_{i,j}$, or $(i,j) \notin T$ and $x \in Q_i$. In this case, by Claim~\ref{eq:U-W-mismatch-zero} and the definition of $W_{\beta}$ we have
\[ \deg_{W_{\beta}}( x ; Q_j ) = d_{i,j}\pi(Q_j), \]
as required for~\ref{itm:graphon-frac-iso-deg-conc}.

Now we shall verify that $W_{\beta}$ satisfies~\ref{itm:graphon-edit-close}. For each $(i,j)\in[M]^2$ define the quantities
\begin{align*}
K_{i,j} & = \int_{(x,y) \in P_i \times P_j} |U(x,y)-U_{\beta}(x,y)| \D\pi_X^2\;\textrm{ and} \\
L_{i,j} & = \int_{(x,y) \in Q_i \times Q_j} |W(x,y)-W_{\beta}(x,y)| \D\pi^2 \;.
\end{align*}
By the definitions of $T'$, $U_{\beta}$ and $W_{\beta}$, we have 
\begin{equation} \label{eq:shave-away-one}
\sum_{(i,j) \in T'} L_{i,j} \le 2\beta\;.
\end{equation}
By Claim~\ref{eq:U-W-mismatch-zero} and the definitions of $U_{\beta}$ and $W_{\beta}$, we have
\begin{equation} \label{eq:quotient-int-upper}
L_{i,j} \le K_{i,j}\textrm{ for all }(i,j) \in [M]^2 \sm T'\;.
\end{equation}
Now by accounting for the editing needed to obtain $U_{\beta}$ from $U$, we obtain
\begin{equation} \label{eq:irregular-shavezero-regsmooth}
\sum_{(i,j) \in [M]^2 \sm T'} K_{i,j} \le \delta + \beta + 4\sqrt{\delta}.
\end{equation}
By definition we have
\[ \Lone{W-W_{\beta}} = \sum_{(i,j)\in[M]^2} L_{i,j}\;, \]
so by combining with equations~\eqref{eq:shave-away-one}--\eqref{eq:irregular-shavezero-regsmooth} we obtain
\[ \Lone{W-W_{\beta}} \le \delta + \beta + 4\sqrt{\delta} + 2\beta \le 4\beta \]
as required.
\end{proof}

Now we prove Lemma~\ref{lem:sampling-outcome}.

\begin{proof}[Proof of Lemma~\ref{lem:sampling-outcome}]
Sample $F\sim\G(m,W)$. Let $\cx = \{X_i\}_{i\in[M]}$ be the partition of $V(F)$ induced (in the sense of Section~\ref{ssec:sampling}) by $\cq$. Let $\cE_1$, $\cE_2$ and $\cE_3$ be the events that~\ref{itm:F-parts-nearly-equal}, \ref{itm:F-frac-iso-deg-nearly-equal} and~\ref{itm:F-random-close} respectively hold. We shall show that $\cE_1 \cap \cE_2 \cap \cE_3$ holds with positive probability, thereby achieving the desired outcome.

We consider the failure probability of $\cE_2$. Let $(i,j)\in[M]^2$ be such that $\mathbf{v}_j\ge \frac{\delta}{M}$. Since the failure probability on the pair $(i,j)$ is trivially zero whenever $d_{i,j} = 0$, we may assume that $d_{i,j}\ge\beta$. Consider $h\in[m]$. Without loss of generality, we shall argue for $h = m$; the argument is entirely analogous for the other choices of $h$. For each $\ell\in[m-1]$ we set $Y_\ell$ to be the indicator function for the event that both $\ell \in X_j$ and $m\ell\in E(F)$ hold; set $Y = \sum_{h=1}^{m-1} Y_h = \deg_F(m;X_j)$. We work in the conditional probability space given $\{m \in X_i\}$. To obtain probability bounds in this conditional probability space, we further condition on $\{x_m=x\}$ for an arbitrary $x\in Q_i$. In the latter conditional probability space, the random variables $(Y_h)_{h\in[m-1]}$ are independent and identically distributed Bernoulli random variables with a common expectation $\mu = \E Y_h$ expressed by~\eqref{eq:W-cell-degree-bound}. Hence, by Theorem~\ref{thm:chernoff-JLR} the probability that $Y \ne (m-1) (1 \pm 6\lambda)d_{i,j}\mathbf{v}_j$ is at most $2\exp\left(-\frac{\lambda^2d_{i,j}\mathbf{v}_jm}4\right)$. It follows that the probability that $Y\ne (m-1) (1 \pm 6\lambda)d_{i,j}\mathbf{v}_j$ is at most $2\exp\left(-\frac{\lambda^2d_{i,j}\mathbf{v}_jm}4\right)$ in the initial conditional probability space. A union bound over $(i,j) \in [M]^2$ and $h\in[m]$ tells us $\cE_2$ does not hold with probability at most $2M^2m\exp\left(-\frac{\lambda^2d_{i,j}\mathbf{v}_jm}4\right) \le 2M^2m\exp\left(-\frac{\lambda^2\beta\mathbf{v}_jm}4\right) < 0.1$.

Now we consider $\cE_1$. For each $i\in[M]$ we have $\pi(Q_i) = \mathbf{v}_i$. Hence, by Theorem~\ref{thm:chernoff-JLR} and a union bound, the event $\cE_1$ fails with probability at most $2M\exp\left(-\frac{\lambda^2m\mathbf{v}_i}3\right)<0.1$. Last, Lemma~\ref{lem:sampling} tells us that the event $\cE_3$ fails with probability at most $\exp\left(-\frac{m}{2 \log_2 m}\right)<0.1$. Hence, $\cE_1 \cap \cE_2 \cap \cE_3$ fails to hold with probability at most 0.3, completing the proof.
\end{proof}

Finally, we prove Lemma~\ref{lem:buffer-precise}.

\begin{proof}[Proof of Lemma~\ref{lem:buffer-precise}]
We construct $(G,\cz)$ from $(F,\cx)$ as follows. We start with $(F,\cx)$. For each $i\in[M]$ we introduce a set $Y_i$ of size 
\begin{equation}\label{eq:howbigY}
|Y_i| = \llbracket (1 +\beta)m\mathbf{v}_i\rrbracket_\Gamma- |X_i|
\end{equation}
and set $Z_i = X_i \sqcup Y_i$. We also introduce a set $Z_0$ of size $|Z_0| = n - \sum_{i\in [M]}|Z_i| \ge 0$.\footnote{The only purpose of $Z_0$, which consists of isolated vertices, is to achieve a given order of the graph.} The first edge-modification resolves~\ref{itm:G-zero} and is easy: we delete all the edges incident to vertices in clusters $Z_i$ with $\mathbf{v}_i<\frac{\delta}M$. Let $F'$ be the graph resulting from these edge deletions.

From now on, we shall work exclusively with pairs in $T=\{(i,j)\in[M]^2:\mathbf{v}_i,\mathbf{v}_j\ge\frac{\delta}M\}$. For each such pair we set 
\begin{equation} \label{eq:targetdeg}
D_{i,j} = d^*_{i,j}|Z_j|\;,
\end{equation}
that is, $D_{i,j}$ is the target quantity on the right-hand side of~\eqref{eq:G-frac-iso-deg-equal}. Note that $D_{i,j}\in \N$. Also, since the matrix $\cd$ is symmetric, so is $\cd^*$ and hence
\begin{equation}\label{eq:ittooktime}
D_{i,j}|Z_i|=d^*_{i,j}\cdot|Z_j|\cdot|Z_i|=d^*_{j,i}\cdot|Z_j|\cdot|Z_i|=D_{j,i}|Z_j|\;.
\end{equation}
Since $\cd^*$ is the $\Gamma$-rounding of $\frac{1 + \alpha}{1+\beta}\cd$ and $\Gamma^{-1}\le\delta^2$, we have
\begin{equation} \label{eq:densitymatrixapprox}
\frac{1+\alpha}{1+\beta}d_{i,j} \le d^*_{i,j} \le (1 + 2\delta)\frac{1+\alpha}{1+\beta}d_{i,j}\;.
\end{equation}
The definition of $\llbracket \cdot \rrbracket_\Gamma$ and the fact that $m\ge \Gamma M/\delta^2$ gives
\begin{equation} \label{eq:vtxpartsize}
(1 - \delta)(1 +\beta)m\mathbf{v}_i \le |Z_i| = \llbracket (1 +\beta)m\mathbf{v}_i\rrbracket_\Gamma \le (1 +\beta)m\mathbf{v}_i\;.
\end{equation}

Let $\tilde{T}:=\{\{i,j\}:(i,j)\in T\cap [M]^{(2)}\}$ be an off-diagonal unordered version of $T$. We add edges in two steps. In the first step, for each $(i,j) \in T$ we add edges between $X_i$ and $Y_j$ to achieve~\eqref{eq:G-frac-iso-deg-equal} for each $v\in X_i$. In the second step, for each $\{i,j\} \in \tilde{T}$ we add edges between $Y_i$ and $Y_j$ to achieve~\eqref{eq:G-frac-iso-deg-equal} for the vertices in $Y_i$ and $Y_j$. Furthermore, for $(i,j) \in T$ with $i=j$ we add edges on $Y_i$ to achieve~\eqref{eq:G-frac-iso-deg-equal} for the vertices in $Y_i$. See Figures~\ref{fig:addingedges} and~\ref{fig:addingedgesinternal} for illustrations of the construction in the cases $i\ne j$ and $i=j$, respectively. Let $G$ be the resultant graph. In particular, \ref{itm:G-zero} still holds in $G$.
\begin{figure}[t]
	\centering
	\includegraphics{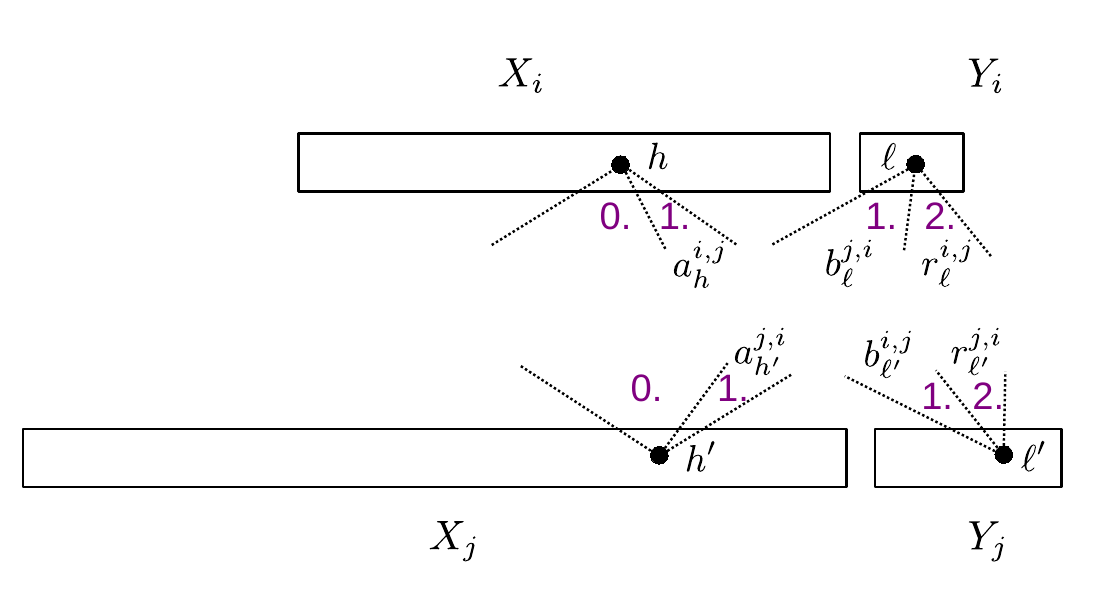}
	\caption{Adding edges in the bipartite graph $(X_i\cup Y_i,X_j\cup Y_j)$ for $i\neq j$. The order of steps is shown in purple (with 0. corresponding to already given edges). The quantities $a_h^{i,j}$, $b^{j,i}_\ell$, $r^{i,j}_\ell$, $a^{j,i}_{h'}$, $b^{i,j}_{\ell'}$ and $r^{j,i}_{\ell'}$ are degree sequences used in the construction.}
	\label{fig:addingedges}
\end{figure}
\begin{figure}[t]
	\centering
	\includegraphics{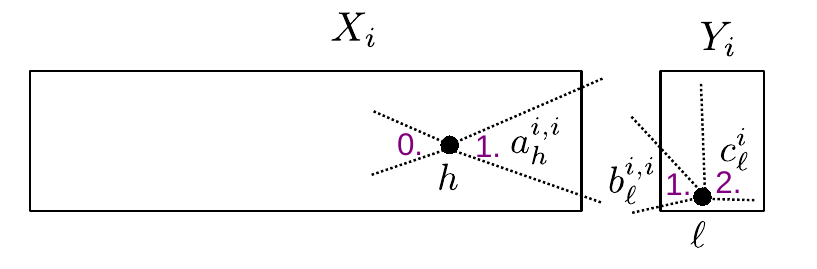}
	\caption{Adding edges in the graph $X_i\cup Y_i$. The order of steps is shown in purple (with 0. corresponding to already given edges). The quantities $a_h^{i,i}$, $b^{i,i}_\ell$, and $c^{i}_\ell$ are degree sequences used in the construction.}
	\label{fig:addingedgesinternal}
\end{figure}

We begin with the first step. For each $h\in X_i$ set $a^{i,j}_h = D_{i,j} - \deg_F(h;X_j)$. Set $A^{i,j} = \sum_{h\in X_i} a^{i,j}_h$. For each $\ell \in Y_j$ set $b^{i,j}_{\ell} \in \left\{\flo{\frac{A^{i,j}}{|Y_j|}},\cei{\frac{A^{i,j}}{|Y_j|}}\right\}$ so that $A^{i,j} = \sum_{\ell\in Y_j}b^{i,j}_{\ell}$. We shall apply the Gale--Ryser theorem (Theorem~\ref{thm:gale-ryser}) to build a bipartite graph with parts $X_i$ and $Y_j$ whose degrees are given by $(a^{i,j}_h)_{h\in X_i}$ and $(b^{i,j}_\ell)_{\ell\in Y_j}$. Fix a non-increasing ordering of $\{a^{i,j}_h\}_{h \in X_i}$ and let $\ca^{i,j}_k$ denote the set of the first $k$ elements of this ordering. We first verify that $a^{i,j}_h\ge0$ for all $h \in X_i$. Indeed, by the definition of $a^{i,j}_h$, \eqref{eq:targetdeg}, \eqref{eq:densitymatrixapprox}, \eqref{eq:vtxpartsize} and~\ref{itm:F-frac-iso-deg-nearly-equal}, for each $h\in X_i$ we have
\begin{equation} \label{eq:predsnyd}
\begin{split}
a^{i,j}_h &\ge (1-\delta)\frac{1+\alpha}{1+\beta}d_{i,j}(1+\beta)\mathbf{v}_jm - (1 + 6\lambda)d_{i,j}\mathbf{v}_jm \\
&\ge d_{i,j}\mathbf{v}_j m(\alpha-7\lambda) \geBy{\eqref{eq:constants2},\eqref{eq:constants4}} 0\;.
\end{split}
\end{equation}
Let us now check the conditions of Theorem~\ref{thm:gale-ryser}. Clearly~\ref{itm:GR1} holds. Furthermore, this implies~\ref{itm:GR2} for $k\ge\max\{b^{i,j}_\ell\}_{\ell\in Y_j}$. Now we focus on the case $k \le \min\{b^{i,j}_\ell\}_{\ell\in Y_j}$; observe that this covers the remaining values of $k$ because the quantities $\{b^{i,j}_{\ell}\}_{\ell\in Y_j}$ differ by at most~$1$. By the definition of $a^{i,j}_h$, \eqref{eq:targetdeg}, \eqref{eq:densitymatrixapprox}, \eqref{eq:vtxpartsize} and~\ref{itm:F-frac-iso-deg-nearly-equal}, for each $h\in X_i$ we have
\begin{equation} \label{eq:snyd}
\begin{split}
a^{i,j}_h &\le (1 + 2\delta)\frac{1 + \alpha}{1+\beta}d_{i,j}(1 + \beta)m\mathbf{v}_j - (1 - 6\lambda)d_{i,j}m\mathbf{v}_j \\
&\le d_{i,j}\mathbf{v}_jm(\alpha+7\lambda) \leByRef{eq:constants4} m \mathbf{v}_j(\beta-7\lambda) \leBy{\ref{itm:F-parts-nearly-equal}} |Y_j|\;.
\end{split}
\end{equation}
Hence, we have $\sum_{h \in \ca^{i,j}_k} a^{i,j}_h \le k|Y_j| \le \sum_{\ell\in Y_j}\min\{k,b^{i,j}_{\ell}\}$. This completes the verification of~\ref{itm:GR2}. Then, by Theorem~\ref{thm:gale-ryser} we have a bipartite graph on $(X_i,Y_j)$ whose edges we add to $F$ to achieve~\eqref{eq:G-frac-iso-deg-equal} for the vertices in $X_i$. Finally, let $F^{\times}$ be the resultant graph after finishing the first step for all $(i,j)\in T$.

Now we consider the second step. Let $\{i,j\}\in\tilde{T}$. For each $h\in Y_i$ set $r^{i,j}_h = D_{i,j} - \deg_{F^{\times}}(h;X_j)$. For each $\ell\in Y_j$ set $r^{j,i}_{\ell} = D_{j,i} - \deg_{F^{\times}}(\ell;X_i)$. We shall apply Theorem~\ref{thm:gale-ryser} to build a bipartite graph with parts $Y_i$ and $Y_j$ with degree sequences $(r^{i,j}_h)_{h\in Y_i}$ and $(r^{j,i}_h)_{h\in Y_j}$. Fix a non-increasing ordering of $\{r^{i,j}_h\}_{h \in Y_i}$ and let $\cb^{i,j}_k$ denote the set of the first $k$ elements of this ordering. We first verify that $r^{i,j}_h,r^{j,i}_{\ell} \ge 0$ for all $h\in Y_i$ and $\ell\in Y_j$. By symmetry, it suffices to prove only the first half of the statement. Every edge between $h \in Y_i$ and $X_j$ was added only in the first step, so we have
\begin{equation} \label{eq:predsnyd2}
\begin{split}
\deg_{F^{\times}}(h;X_j) &= b^{j,i}_h \le \cei{\frac{A^{j,i}}{|Y_i|}}\leByRef{eq:snyd} 1+ \frac{d_{j,i}\mathbf{v}_i m(\alpha+7\lambda)\cdot |X_j|}{\llbracket (1 +\beta)m\mathbf{v}_i\rrbracket_\Gamma- |X_i|}\\
&\leBy{\ref{itm:F-parts-nearly-equal}} \frac{d_{j,i}(\alpha+9\lambda)\mathbf{v}_jm}{\beta-6\lambda-\delta}
\leByRef{eq:constants4} d_{i,j}\mathbf{v}_jm \le D_{i,j}\;.
\end{split}
\end{equation}
This shows that $r^{i,j}_h\ge 0$. Now we verify the conditions of Theorem~\ref{thm:gale-ryser}, starting with~\ref{itm:GR1}. We have
\begin{align*}
\sum_{h \in Y_i}r^{i,j}_h & = D_{i,j}|Y_i| - \sum_{h \in Y_i}b^{j,i}_h = D_{i,j}|Y_i| - \sum_{\ell \in X_j}a^{j,i}_{\ell} \\
& = D_{i,j}|Y_i| - \sum_{\ell\in X_j}\big(D_{j,i}-\deg_{F}(\ell;X_i)\big)\\
& = D_{i,j}|Y_i|-D_{j,i}|X_j| + e_{F}(X_i,X_j)\;.
\end{align*}
By symmetry, we also have $\sum_{h \in Y_j}r^{j,i}_h=D_{j,i}|Y_j|-D_{i,j}|X_i| + e_{F}(X_i,X_j)$. Now, crucially, recall that by~\eqref{eq:ittooktime} we have $D_{i,j}(|X_i|+|Y_i|)=D_{j,i}(|X_j|+|Y_j|)$. Hence, we obtain $\sum_{h \in Y_i}r^{i,j}_h=\sum_{h \in Y_j}r^{j,i}_h$ as needed for~\ref{itm:GR1}. Furthermore, this implies~\ref{itm:GR2} for $k \ge \max\{r^{j,i}_\ell\}_{\ell \in Y_j}$. Now we focus on the case $k \le \min\{r^{j,i}_\ell\}_{\ell \in Y_j}$. Observe that this covers the remaining values of $k$ because we have $r^{j,i}_{\ell}=D_{j,i}-b^{i,j}_{\ell}$ and the fact that the quantities $\{b^{i,j}_{\ell}\}_{\ell \in Y_j}$ differ by at most~$1$ means that the quantities $\{r^{j,i}_{\ell}\}_{\ell \in Y_j}$ also differ by at most~$1$. For each $\ell\in Y_i$ we have
\begin{align*}
\deg_{F^{\times}}(\ell;X_j)&=b^{j,i}_h\ge\flo{\frac{A^{j,i}}{|Y_i|}} \ge \frac{1}{|Y_i|}\sum_{h \in X_j} a^{j,i}_h -1 \\
& \geByRef{eq:predsnyd} \frac{d_{j,i}\mathbf{v}_im(\alpha-7\lambda)\cdot|X_j|}{\llbracket(1+\beta)m\mathbf{v}_i\rrbracket_\Gamma-|X_i|} - 1\\
& \geBy{\ref{itm:F-parts-nearly-equal}} \frac{d_{j,i}(\alpha-9\lambda)\mathbf{v}_jm}{\beta + 6\lambda} \\
& \geByRef{eq:constants4} [d_{i,j}(1+\alpha+2\lambda) - (\beta-7\lambda)]\mathbf{v}_jm \ge D_{i,j} - |Y_j|\;,
\end{align*}
where the fact that $d_{i,j} \le 1-\beta$ is used in the penultimate inequality and the definition of $D_{i,j}$ is utilized in the final inequality. Hence, we have $\sum_{h \in \cb^{i,j}_k} r^{i,j}_h \le k|Y_j| \le \sum_{\ell\in Y_j}\min\{k,r^{j,i}_{\ell}\}$. This completes the verification of~\ref{itm:GR2}. Then, by Theorem~\ref{thm:gale-ryser} we have a bipartite graph on $(Y_i,Y_j)$ whose edges we add to $F^{\times}$ to achieve~\eqref{eq:G-frac-iso-deg-equal} for the vertices in $Y_i$ and $Y_j$.

Now we handle the second step for the case $i=j$. We need to achieve~\eqref{eq:G-frac-iso-deg-equal} in $Y_i$. To this end, we shall apply the Erd\H{o}s--Gallai Theorem (Theorem~\ref{thm:erdos-gallai}) to build a graph on $Y_i$ with degree sequence $\{c^i_\ell\}_{\ell \in Y_i}$ given by $c^i_\ell = D_{i,i} - b^{i,i}_\ell$ for each $\ell \in Y_i$. Observe that for $\Delta = \max\{c^i_\ell\}_{\ell\in Y_\ell} \in \N$ and for all $\ell\in Y_\ell$ we have $c^i_\ell\in\{\Delta-1,\Delta\}$. Fix a non-increasing ordering of $\{c^i_\ell\}_{\ell \in Y_i}$ and let $\cb^{i,i}_k$ denote the set of the first $k$ elements of this ordering. We first verify that $0 \le c^i_\ell \le |Y_i|-1$ for all $\ell \in Y_i$. Noting that every edge between $\ell \in Y_i$ and $X_i$ was added only in the first step and repeating~\eqref{eq:predsnyd2} mutatis mutandis, for each $\ell\in Y_i$ we have
\begin{equation}\label{eq:predsnyd3}
\deg_{F^{\times}}(\ell;X_i) = b^{i,i}_h \le \frac{d_{i,j}(\alpha+8\lambda)|X_j|}{\beta-6\lambda-\delta} \le D_{i,i}\;.
\end{equation}
This shows that $c^i_\ell \ge 0$. Furthermore, we have
\begin{align*}
\deg_{F^{\times}}(\ell;X_i)&=b^{i,i}_h\ge\flo{\frac{A^{i,i}}{|Y_i|}} \geByRef{eq:predsnyd} \frac{d_{i,i}\mathbf{v}_im(\alpha-7\lambda)\cdot|X_i|}{\llbracket(1+\beta)m\mathbf{v}_i\rrbracket_\Gamma-|X_i|} - 1\\
&\geBy{\ref{itm:F-parts-nearly-equal}} \frac{d_{i,j}(\alpha-8\lambda)|X_j|}{\beta + 6\lambda} \\
&\geByRef{eq:constants4} [d_{i,i}(1+\alpha+2\lambda) - (\beta-7\lambda)]\mathbf{v}_im  \ge D_{i,i} - |Y_i| + 2\;,
\end{align*}
where the fact that $d_{i,i} \le 1-\beta$ is used in the penultimate inequality and the definition of $D_{i,i}$ is utilized in the final inequality. This shows that $c^i_\ell \le |Y_i|-2$.

Let us now check the conditions of Theorem~\ref{thm:erdos-gallai}. First, we verify~\ref{itm:EG1}, that is, the quantity $\sum_{\ell \in Y_i} c^i_\ell$ is even. Since $\sum_{h \in X_i}\deg_F(h;X_i)=2e_F(X_i)$ and $D_{i,i}$ are both even, the quantity
\begin{align*}
\sum_{\ell \in Y_i} c^i_\ell = D_{i,i}|Y_i|-\sum_{\ell \in Y_i}b^{i,i}_\ell & = D_{i,i}|Y_i|-\sum_{h \in X_i}a^{i,i}_h \\
& = D_{i,i}|Y_i| - D_{i,i}|X_i| + \sum_{h \in X_i} \deg_F(h;X_i)
\end{align*}
is even as needed. Now we turn to~\ref{itm:EG2}. We shall consider three different cases. First, we consider the case $k > \Delta$. In this case, we have 
\[
\sum_{\ell \in \cb^{i,i}_k} c^i_\ell \le k\Delta \le k(k-1)\;,
\]
which implies~\ref{itm:EG2} as required. Second, we have the case $k < \Delta$. In this case, we have 
\[
\sum_{\ell \in \cb^{i,i}_k} c^i_\ell \le k(|Y_i|-1) \le k(k-1) + (|Y_i|-k)\min\{k,\Delta - 1\}\;,
\]
which implies~\ref{itm:EG2} as required. Finally, we have the case $k = \Delta \le |Y_i| - 2$. Observe that in this case we have
\[
\sum_{\ell \in \cb^{i,i}_k} c^i_\ell \le k(k-1) + \min\{k,c_{k+1}\} + \min\{k,c_{k+2}\}\;,
\]
except in the case where $c_k > c_{k+1} = c_{k+2} = 0$. Observe that in the exceptional case we must have $k = \max\{c^i_\ell\} = 1$ and hence $\sum_{\ell \in Y_i} c^i_\ell = 1$, so~\ref{itm:EG1} implies that it never holds. Hence, the inequality above holds, implying~\ref{itm:EG2} as required. Then, by Theorem~\ref{thm:erdos-gallai} we have a graph on $Y_i$ whose edges we add to $F^{\times}$ to achieve~\eqref{eq:G-frac-iso-deg-equal} in $Y_i$.

To complete the proof, we shall check~\ref{itm:G-F-close}. Since $1+\beta \le \frac{n}{m} \le 1+2\beta$ and $G[V(F)]=F'$, by Lemma~\ref{lem:zoom} we have $\cutm(G,F')\le2\left(1-\frac{m}{n}\right)\le\frac{4\beta}{1+2\beta}\le4\beta$. Thus, it remains to prove that $\cutm(F,F')\le\beta$. To do this, we recall that the only edge edits wholly contained in $V(F)$ are the deletions of edges incident to clusters $X_i$ with $\mathbf{v}_i<\frac{\delta}{M}$. Each such cluster has at most $|X_i|=(1\pm\lambda)\mathbf{v}_i m<\frac{2\delta n}M$ vertices, so at most $M\times \frac{2\delta n}M\times n<\beta n^2$ edges are deleted. This completes the proof.
\end{proof}

\section{Regular case}\label{sec:regular}

As mentioned in the introduction, the proof of Theorem~\ref{thm:main-approx} given in Section~\ref{sec:proof} guarantees regularity of the approximating graphs in the case that the graphons are regular, which leads us to Theorem~\ref{thm:regular-approx}. Recall that a graphon is regular with degree $d$ if and only if it is fractionally isomorphic to the constant-$d$ graphon.

We shall walk through each of the three steps (Lemmas~\ref{lem:graphon-frac-iso-approx}--\ref{lem:buffer-precise}) of the proof of Theorem~\ref{thm:main-approx} and see that they indeed yield Theorem~\ref{thm:regular-approx}. Since the cases $d\in\{0,1\}$ are trivial, we may assume that $d\in(0,1)$. In particular, we may choose the constant $\beta>0$ to satisfy $\beta\ll \min\{d,1-d\}$. We may further assume that $C_d\in\cu$. 
Let us first deal with Lemma~\ref{lem:graphon-frac-iso-approx}. That is, we pick a representative of $\cu$ to play the role of $W'$ in Lemma~\ref{lem:graphon-frac-iso-approx}; here we shall make the obvious choice of $W'=C_d$. In particular, $W'/\cc(W')$ is a constant-$d$ graphon defined on a singleton probability space $\Omega/\cc(W')$. Hence, $M=1$ and the $1\times 1$-matrix $\cd$ contains a single entry of value $d$. Furthermore, each $W\in\cu$ satisfies~\ref{itm:graphon-frac-iso-deg-conc} automatically. Now Lemma~\ref{lem:sampling-outcome} is an intermediate step which yields for each $W\in\cu$ a graph $F_W$ on $m$ vertices such that, by~\ref{itm:F-frac-iso-deg-nearly-equal}, each vertex $x\in V(F_W)$ satisfies
\begin{equation}\label{eq:Iwhst1}
\deg_{F_W}(x) = (1 \pm \eps)dm\;.
\end{equation}
We pick $m\approx\frac{n}{1+\beta}$. Finally, we deal with Lemma~\ref{lem:buffer-precise}. Here we create graphs $G_W$ by adding some vertices and edges to $F_W$. More precisely, we fix an even integer $D\approx(1+\beta)dm$ and add a little set $Y_1$ to the set $X_1=V(F_W)$. Then, we add edges onto both the pair $(X_1,Y_1)$ (one application of the Gale--Ryser Theorem) and inside the set $Y_1$ (one application of the Erd\H{o}s--Gallai Theorem) so that the graph becomes regular with degree $D$. Observe that, in the proof of Lemma~\ref{lem:buffer-precise}, the set $Y_1$ is created in~\eqref{eq:howbigY} in a way that $|X_1\cup Y_1|$ is divisible by $\Gamma$. However, this divisibility condition can be dropped in our setting as its only purpose\footnote{as explained after the statement of Lemma~\ref{lem:sampling-outcome}} was to allow for biregularity between $X_1\cup Y_1$ and other clusters $X_i\cup Y_i$, of which there are none in the current setting. Therefore, we are unconstrained by divisibility considerations and have enough flexibility to choose $|Y_1|$ so as to achieve $|X_1\cup Y_1| = n$. Checking the conditions of the Gale--Ryser theorem and of the Erd\H{o}s--Gallai Theorem works the same way as before.

In addition to the steps mentioned above, the proof of Theorem~\ref{thm:main-approx} involves two further steps in which edges incident to small clusters are deleted and a set $Z_0$ of isolated vertices is added. We remark that the former step does not apply in our setting. Furthermore, the addition of the set $Z_0$ serves the sole purpose of ensuring that the graph $G$ has a specific number of vertices; this is achieved in the current setting by having more flexibility in the size of the set $Y_1$.

\section*{Acknowledgments}
We thank two anonymous referees  for their comments.

\bibliographystyle{plain}
\bibliography{DSG}

\end{document}